\documentclass[12pt]{article}
\usepackage{latexsym, amsmath, amsfonts, amscd, amssymb, verbatim}
\begin{document}
\pagestyle{myheadings}

\newtheorem{Theorem}{Theorem}[section]
\newtheorem{Proposition}[Theorem]{Proposition}
\newtheorem{Remark}[Theorem]{Remark}
\newtheorem{Lemma}[Theorem]{Lemma}
\newtheorem{Corollary}[Theorem]{Corollary}
\newtheorem{Definition}[Theorem]{Definition}
\newtheorem{Example}[Theorem]{Example}
\renewcommand{\theequation}{\thesection.\arabic{equation}}
\normalsize

\setcounter{equation}{0}

\title{\bf A Representation of Generalized Convex Polyhedra and Applications\footnote{This work was supported by National Foundation for Science $\&$ Technology Development (Vietnam) under the grant No.~101.01-2014.37. The second author thanks the Vietnam Institute for Advanced Study in Mathematics for supporting his 6-month stay at the Institute in 2015.}}

\medskip
\author{Nguyen Ngoc Luan\footnote{Department of
Mathematics and Informatics, Hanoi National University of Education, 136 Xuan Thuy, Hanoi, Vietnam; email: luannn@hnue.edu.vn.}\ \, and\ \,
Nguyen Dong Yen\footnote{Institute of Mathematics, Vietnam Academy of
Science and Technology, 18 Hoang Quoc Viet, Hanoi 10307, Vietnam;
email: ndyen@math.ac.vn.}}\maketitle

\medskip
\begin{quote}
\noindent {\bf Abstract.} It is well known that finite-dimensional polyhedral convex sets can be generated by finitely many points and finitely many directions. Representation formulas in this spirit are obtained for convex polyhedra and generalized convex polyhedra in locally convex Hausdorff topological vector spaces. Our results develop those of X.~Y.~Zheng (Set-Valued Anal., Vol. 17, 2009, 389--408), which were established in a Banach space setting. Applications of the representation formulas to proving solution existence theorems for generalized linear programming problems and generalized linear vector optimization problems are shown.



\end{quote}

\newpage
\section{Introduction}
The intersection of a finite number of closed half-spaces of a finite-dimensional Euclidean space is called a \textit{polyhedral convex set} (a \textit{convex polyhedron} in brief). By convention, the intersection of an empty family of closed half-spaces is the whole space. Therefore, emptyset and the whole spaces are two special polyhedra.  Due to \cite[Theorem 19.1]{Rock_book_1970}, for every given convex polyhedron one can find a finite number of points and a finite number of directions such that the polyhedron can be represented as the sum of the convex hull of those points and the convex cone generated by those directions. The converse is also true. This celebrated result is attributed \cite[p. 427]{Rock_book_1970} primarily to Minkowski \cite{Minkowski_1910} and Weyl \cite{Weyl_1935, Weyl_1953}. By using the result, it is easy to derive fundamental solution existence theorems in linear progamming. Note that the just cited representation formula for finite-dimensional polyhedral convex sets has many other applications in mathematics. As an example, one can refer to the elegant proofs of the necessary and sufficient second-oder conditions for a local solution and for a locally unique solution in quadratic programming, which were given by Contesse \cite{Contesse_1980} in 1980; see \cite[pp.~50--63]{Lee_Tam_Yen_2005} for details.

\medskip
According to Bonnans and Shapiro \cite[Definition 2.195]{Bonnans_Shapiro_2000}, a subset of locally  convex Hausdorff topological vector space is called a \textit{generalized polyhedral convex set} (or a \textit{generalized convex polyhedron}) if it is the intersection of finitely many closed half-spaces and a closed affine subspace of that topological vector space. If the affine subspace can be chosen as the whole space, the generalized polyhedral convex set is said to be a \textit{polyhedral convex set} (or a \textit{convex polyhedron}). The theories of generalized linear programming and quadratic programming in \cite[Sections 2.5.7 and 3.4.3]{Bonnans_Shapiro_2000} are based on this concept of generalized convex polyhedron. 

\medskip
In 2009, using a result related to the Banach open mapping theorem (see, e.g., \cite[Theorem 5.20]{Rudin_1991}), Zheng \cite[Corollary ~2.1]{Zheng_2009} has clarified the relationships between convex polyhedra in Banach spaces and the finite-dimensional convex polyhedra. 

\medskip   
It is well known that any infinite-dimensional normed space equipped with the {\it weak topology} is not metrizable, but it is a locally convex Hausdorff topological vector space. Similarly, the dual space of any infinite-dimensional normed space equipped with the {\it weak$^*$ topology} is not metrizable, but it is a locally convex Hausdorff topological vector space. Actually, the just mentioned two models provide us with the most typical examples of locally convex Hausdorff topological vector spaces, whose topologies cannot be given by norms. It is clear that Zheng's results in \cite{Zheng_2009} cannot be used neither for a infinite-dimensional normed space equipped with the weak topology, nor for the dual space of any infinite-dimensional normed space equipped with the weak$^*$ topology.

\medskip
The aim of our paper is twofold: to find an analogue of the above-mentioned representation of finite-dimensional convex polyhedra via finite families of points and directions for convex polyhedra in locally convex Hausdorff topological vector spaces, and to apply the obtained results to proving solution existence theorems for infinite-dimensional linear programming problems and linear vector optimization problems. Among other things, we will show that the result of  Zheng \cite[Corollary ~2.1]{Zheng_2009} is valid for convex polyhedra in locally convex Hausdorff topological vector spaces. 

\medskip
The organization of the present paper is as follows. In Section 2, we obtain representation formulas for generalized convex polyhedra. Section 3 is devoted to solution existence of generalized linear programs. Solution existence of generalized linear vector optimization problems is studied in Section 4.

\section{Representation Formulas for Generalized  Convex Polyhedra}
Let $X$ be a {\it locally convex Hausdorff topological vector space} with the dual space denoted by $X^*$. For any $x^* \in X^*$ and $x \in X$, $\langle x^*, x \rangle$ indicates the value of $x^*$ at $x$. 
\begin{Definition}{\rm (See \cite[p.~133]{Bonnans_Shapiro_2000})
A subset $C \subset X$ is said to be a \textit{generalized polyhedral convex set} (a \textit{generalized convex polyhedron} for short) if there exist $x^*_i \in X^*$, $\alpha_i \in \mathbb R$, $i=1,2,\dots,p$, and a closed affine subspace $L \subset X$, such that 
\begin{equation*}
C=\left\{ x \in X \mid x \in L, \; \langle x^*_i, x \rangle \leq \alpha_i,\ \;  i=1,\dots,p\right\}.
\end{equation*} If $C$ admits the last representation for $L=X$ and for some $x^*_i \in X^*$, $\alpha_i \in \mathbb R$, $i=1,2,\dots,p$, then it is called a \textit{polyhedral convex set} (or a \textit{convex polyhedron}).}
\end{Definition}

The following classical result shows that, for any convex polyhedron in $\mathbb R^n$, one can find a finite number of points and a finite number of directions such that the polyhedron can be represented as the sum of the convex hull of those points and the convex cone generated by those directions. The converse is also true.
\begin{Theorem}{\rm  (\cite[Theorem~19.1]{Rock_book_1970})}\label{Rock1970} For any nonempty set $C \subset \mathbb R^n$, the following two properties are equivalent:\\
\indent {\rm (a)} $C$ is a convex polyhedron;\\
\indent {\rm (b)} $C$ is finitely generated, i.e., $C$ can be represented as
\begin{equation}\label{rep_1}
\begin{aligned}
C=\Bigg\{ \sum\limits_{i=1}^k \lambda_i u_i + \sum\limits_{j=1}^\ell \mu_j v_j \mid & \ \lambda_i \geq 0, \ \forall i=1,\dots,k,    \\ 
&\sum\limits_{i=1}^k \lambda_i=1,\ \, \mu_j \geq 0,\ \forall j=1,\dots,\ell   \Bigg\},&& 
\end{aligned}
\end{equation}
where $u_i \in \mathbb R^n$, $i=1, \dots, k$, and $v_j \in \mathbb{R}^n$, $j=1,\dots,\ell$.
\end{Theorem}

From \eqref{rep_1} it follows that $u_i \in C$ for $i=1, \dots, k$.
A natural question arises: \textit{Is there any analogue of the representation \eqref{rep_1} for convex polyhedra in locally convex Hausdorff topological vector spaces, or not?} In order to give an answer in the affirmative to this question, we will need several results from functional analysis. In what follows, $X$ is a locally convex Hausdorff topological vector space.

\begin{Lemma}{\rm (Closedness of the sum two linear subspaces; see  \label{Rudin1991} \cite[Theorem 1.42]{Rudin_1991})} 
Suppose $X_0$ and $X_1$ are linear subspaces of $X$, $X_0$ is closed, and $X_1$ has finite dimension. Then $X_0+X_1$ is closed.   
\end{Lemma}

\begin{Lemma}{\rm (The Hahn-Banach extension theorem; see \label{Hann_Banach} \cite[Theorem 3.6]{Rudin_1991})} 
If $x^*$ is a continuous linear functional on a linear subspace $M$ of $X$, then there exists $\widetilde{x}^* \in X^*$ such that $\langle \widetilde{x}^*, x \rangle = \langle x^*, x \rangle $ for all $x \in M$.   
\end{Lemma}

The forthcoming lemma follows from a theorem in \cite{Rudin_1991}. A proof is provided here for the sake of clarity of our presentation. 

\begin{Lemma}{\rm \label{cont_proj}} 
	If $Y$ and $Z$ are Hausdorff finite-dimensional topological vector spaces of dimension $n$ and if $g: Y \rightarrow Z$ is a linear bijective mapping, then $g$ is a homeomorphism.
\end{Lemma}
{\bf Proof.} Let $\{e_1, e_2, \dots, e_n\}$ be a basis of the Euclidean space ${\mathbb R}^n$, which is equipped with the natural topology. Let $\{v_1, v_2, \dots, v_n\}$ be a basis of $Y$. Setting $w_i=g(v_i)$ for $i=1,\dots,n$, we see that $\{w_1, w_2, \dots, w_n\}$ is a basis of $Z$. Clearly, there is an unique linear bijection $\Phi: {\mathbb R}^n \rightarrow Y$ satisfying the conditions $\Phi(e_i)=v_i$ for all $i$. Similarly, there is an unique linear bijection $\Psi: {\mathbb R}^n \rightarrow Z$ with $\Phi(e_i)=w_i$ for all $i$. By \cite[Theorem 1.21(a)]{Rudin_1991}, $\Phi$ and $\Psi$ are homeomorphisms. (Note that the quoted result was obtained for ${\mathbb C}^n$ and topological vector spaces over the complex field $\mathbb C$. Nevertheless, the method of proof is valid for the case of ${\mathbb R}^n$ and topological vector spaces over $\mathbb R$.)  Since $g=\Psi \circ \Phi^{-1}$ and $g^{-1}=\Phi \circ \Psi^{-1}$ by our construction, it follows that both $g$ and $g^{-1}$ are continuous mappings. $\hfill\Box$
 
\medskip 
 We are now in a position to extend Corollary 2.1 from the paper of Zheng \cite{Zheng_2009}, which was given in a normed spaces setting, to the case of convex polyhedra in locally convex Hausdorff topological vector spaces.
\begin{Proposition}\label{Zheng_ext}
A nonempty subset $D \subset X$ is a convex polyhedron if only if there exist closed linear subspaces $X_0$, $X_1$ of $X$ and a convex polyhedron $D_1 \subset X_1$ such that
\begin{equation}\label{decomp_for_X}
X=X_0+X_1, \quad X_0 \cap X_1=\{0\},\quad   {\rm{dim}}\,X_1 < +\infty,
\end{equation}
 and 
\begin{equation}\label{decomp_for_D}
D=D_1+X_0.
\end{equation}
\end{Proposition}
{\bf Proof.} {\it Necessity:} If $D$ is a convex polyhedron, then there exist $x^*_i \in X^*$, $\alpha_i \in \mathbb R$, $i=1,\dots,p$, such that
\begin{equation*}
D=\left\{ x \in X \mid \langle x^*_i, x \rangle \leq \alpha_i, \; i=1,\dots,p\right\}.
\end{equation*}
Let 
\begin{equation*}
X_0:=\left\{ x \in X \mid \langle x^*_i, x \rangle =0, \; i=1,\dots,p\right\}.
\end{equation*}
Because $X_0$ is a closed linear subspace of finite codimension, one can find a finite-dimensional linear subspace $X_1$ of $X$, such that $X=X_0 + X_1$ and $X_0 \cap X_1=\{0\} $.  By \cite[Theorem 1.21(b)]{Rudin_1991}, $X_1$ is closed. Clearly, $$D_1:=\left\{ x \in X_1 \mid \langle x^*_i, x \rangle \leq \alpha_i, \; i=1,\dots,p\right\}$$ is a convex polyhedron in $X_1$. It is easy to verify that $D_1 + X_0 \subset D$. The reverse inclusion is also true. Indeed, for each $x \in D$ there exist $x_0 \in X_0$ and $x_1 \in X_1$ satisfying $x=x_0+x_1$. 
Since
\begin{equation*}
\langle x^*_i, x_1 \rangle = \langle x^*_i, x \rangle - \langle x^*_i, x_0 \rangle = \langle x^*_i, x \rangle \leq \alpha_i
\end{equation*}  
 for all $i=1,\dots,p$, it follows that $x_1 \in D_1$; hence $x=x_1+x_0 \in D_1+X_0$. We have thus proved that $D=D_1+X_0$.

{\it Sufficiency:}  Let  $X_0$, $X_1$ be closed subspaces of $X$ satisfying the conditions in \eqref{decomp_for_X}. Let $D_1 \subset X_1$ be a convex polyhedron in $X_1$ and let $D$ be defined by \eqref{decomp_for_D}. Select  $u^*_j \in X_1^*$ and $\beta_j \in \mathbb R$, $j=1,\dots,m$, such that
\begin{equation*}
D_1=\left\{ u \in X_1 \mid \langle u^*_j, u \rangle \leq \beta_j, \; j=1,\dots,m\right\}.
\end{equation*}
Let $\pi_0: X \rightarrow X/X_0$, $x \mapsto x+X_0$ for all $x \in X$, be the canonical projection from $X$ on the quotient space $X/X_0$. It is clear that the operator $\Phi_0: X/X_0 \rightarrow X_1$, $x_1 +X_0 \mapsto x_1$ for all $x_1 \in X_1$, is a linear bijective mapping. On one hand, by \cite[Theorem 1.41({\it a})]{Rudin_1991}, $\pi_0$ is a linear continuous mapping. On the other hand, $\Phi_0$ is a homeomorphism by Lemma \ref{cont_proj}. So, the operator $\pi:=\Phi_0\circ \pi_0: X \rightarrow X_1$ is linear and continuous. Put $x^*_j=u^*_j \circ \pi$, $j=1,\dots,m$. Take any $x=x_1+x_0$ with $x_1 \in D_1$ and $x_0 \in X_0$. It clear that
\begin{equation*}
\langle x^*_j, x \rangle = \langle u^*_j \circ \pi , x \rangle =\langle u^*_j,  \pi(x) \rangle =  \langle u^*_j, x_1 \rangle  \leq \beta_j
\end{equation*}  
for all $j=1,\dots,m$. Conversely, take any $x \in X$ satisfying $\langle x^*_j, x \rangle \leq \beta_j$ for all $j=1,\dots,m$. Let  $x_0 \in X_0$ and $x_1 \in X_1$ be such that $x=x_0+x_1$. Since $$\beta_j \geq \langle x^*_j, x_0+x_1 \rangle  = \langle u^*_j \circ \Phi_0\circ \pi_0, x_0+x_1 \rangle =   \langle u^*_j, x_1 \rangle $$ for all $j=1,\dots,m$, we see that $x_1 \in D_1$.  Hence $x \in D_1 + X_0$. It follows that 
$D_1+X_0=\left\{x \in X \mid \langle x^*_j, x \rangle \leq \beta_j, \; j=1,\dots,m\right\}$. Therefore $D=D_1+X_0$ is a convex polyhedron in $X$. $\hfill\Box$

\medskip
The main result of this section is formulated as follows.

\begin{Theorem}\label{main_result}
A nonempty subset $D \subset X$ is a generalized convex polyhedron if and only if there exist $u_1, \dots, u_k \in X$, $v_1, \dots, v_{\ell} \in X$, and a closed linear subspace  $X_0 \subset X$ such that
\begin{equation}\label{rep_2}
\begin{aligned}
D=\Bigg\{ \sum\limits_{i=1}^k \lambda_i u_i + \sum\limits_{j=1}^\ell \mu_j v_j \mid & \lambda_i \geq 0, \ \forall i=1,\dots,k,    \\ 
&\sum\limits_{i=1}^k \lambda_i=1,\ \, \mu_j \geq 0,\ \forall j=1,\dots,\ell   \Bigg\}+X_0.&& 
\end{aligned}
\end{equation}
\end{Theorem}
{\bf Proof.} {\it Necessity:} Suppose that $D$ is a generalized convex polyhedron. Then we have 
\begin{equation*}
D=\left\{x \in X \mid x \in L, \ \langle x^*_i, x \rangle \leq \alpha_i, \ i=1,2,\dots,p\right\},
\end{equation*}
where $L \subset X$ is a closed affine subspace, $x^*_i \in X^*$ and $\alpha_i \in \mathbb R$ for $i=1,\dots,p$. Select a locally convex Hausdorff topological vector
space $Y$, a continuous linear mapping $A: X \rightarrow Y$, and a point $y \in Y$ such that $L=\{x \in X \mid Ax=y\}$. Fix an element $x_0 \in D$ and set $D_0=D-x_0$. It is easy to verify that
\begin{equation*}
D_0=\left\{ u \in X \mid Au=0, \ \, \langle x^*_i, u \rangle \leq \alpha_i-\langle x^*_i, x_0 \rangle, \ i=1,\dots,p\right\}.
\end{equation*}
As $D_0$ is a convex polyhedron in ${\rm ker}A:=\{u \in X \mid Au=0\}$, by Proposition \ref{Zheng_ext} we can find closed linear subspaces $X_{0, A}$ and $X_{1, A}$ of ${\rm{ker}}A$ and a convex polyhedron $D_{1, A} \subset X_{1, A}$ such that
\begin{equation*}
{\rm{ker}}A=X_{0, A}+X_{1, A}, \quad X_{1, A} \cap X_{0, A}=\{0\},\quad  {\rm{dim}}X_{1, A} < +\infty,
\end{equation*}
and 
\begin{equation*}
D_0=D_{1, A}+X_{0, A}.
\end{equation*}
Because $X_{1,A} \subset {\rm{ker}}A$ is closed and ${\rm{ker}}A$ is a closed linear subspace of $X$, $X_{1,A}$ is a closed linear subspace of $X$. Since $D_{1, A}$ is a convex polyhedron of the finite-dimensional space $X_{1, A}$, invoking Theorem~\ref{Rock1970} we can represent $D_{1, A}$ as
\begin{equation*}
\begin{aligned}
D_{1, A}=\Bigg\{ \sum\limits_{i=1}^k \lambda_i u_i + \sum\limits_{j=1}^\ell \mu_j v_j \mid & \lambda_i \geq 0, \ \forall i=1,\dots,k,    \\ 
&\sum\limits_{i=1}^k \lambda_i=1,\ \, \mu_j \geq 0,\ \forall j=1,\dots,\ell   \Bigg\},&& 
\end{aligned}
\end{equation*}
where $u_i \in D_{1, A}$ for $i=1,\dots,k$ and $v_j \in X_{1, A}$ for $j=1,\dots,\ell$.
Therefore 
\begin{equation*}
\begin{aligned}
D=\Bigg\{ \sum\limits_{i=1}^k \lambda_i \left({u}_i+x_0\right) + \sum\limits_{j=1}^\ell \mu_j {v}_j \mid &\ \lambda_i \geq 0, \ \forall i=1,\dots,k,    \\ 
&\sum\limits_{i=1}^k \lambda_i=1,\ \, \mu_j \geq 0,\ \forall j=1,\dots,\ell   \Bigg\}+{X}_{0, A}.&& 
\end{aligned}
\end{equation*}
We have thus found a representation of the form \eqref{rep_2} for $D$.

{\it Sufficiency:} Suppose that $D$ is of the form \eqref{rep_2}. 
Let $$X_1={\rm{span}}\{u_1,\dots,u_k,v_1,\dots,v_{\ell}\}$$ be the linear subspace generated by the vectors $u_1,\dots,u_k,v_1,\dots,v_{\ell}$. Put
\begin{equation*}
\begin{aligned}
D_1:=\Bigg\{ \sum\limits_{i=1}^k \lambda_i u_i + \sum\limits_{j=1}^\ell \mu_j v_j \mid & \lambda_i \geq 0, \ \forall i=1,\dots,k,    \\ 
&\sum\limits_{i=1}^k \lambda_i=1, \mu_j \geq 0, \ \forall j=1,\dots,\ell   \Bigg\}.&& 
\end{aligned}
\end{equation*}
By Lemma \ref{Rudin1991}, $W:=X_1 + X_0$ is a closed linear subspace of $X$. 
Because $X_0$ is a closed subspace of finite codimension of $W$, one can find a finite-dimensional linear subspace $W_1 \subset W$, such that $W=X_0 + W_1$ and $X_0 \cap W_1=\{0\} $. Consider the continuous linear mapping $\pi: W \rightarrow W_1$ be defined by $\pi(x)=w_1$, where $x=x_0+w_1$ with $\left(w_1, x_0\right) \in W_1 \times X_0$.
We have
\begin{equation*}
\begin{aligned}
\pi\left(D_1\right)=\Bigg\{ \sum\limits_{i=1}^k \lambda_i \pi\left(u_i \right) + \sum\limits_{j=1}^\ell \mu_j \pi\left(v_j\right) \mid & \lambda_i \geq 0, \ \forall i=1,\dots,k,    \\ 
&\sum\limits_{i=1}^k \lambda_i=1,\ \,\mu_j \geq 0,\ \forall j=1,\dots,\ell   \Bigg\}.&& 
\end{aligned}
\end{equation*}
By Lemma \ref{Rock1970}, $\pi(D_1)$ is a convex polyhedron of $W_1$. We have  $D_1+X_0=\pi(D_1)+X_0$. Indeed, if $x=x_1+x_0$ where $x_1 \in D_1$ and $x_0 \in X_0$, then $x_{1,0}:=x_1-\pi(x_1) \in X_0$. So $x=\pi(x_1)+x_{1,0}+x_0 $ belongs to $\pi(D_1)+X_0$. Conversely, for any $x=\pi(z_1)+x_0$ with $z_1 \in D_1$ and $x_0 \in X_0$, we have 
\begin{equation*}
x=z_1 + \pi(z_1) - z_1 + x_0=z_1 + \left(x_0 - \left(z_1-\pi(z_1)\right)\right) \in D_1 + X_0.
\end{equation*}
Since $D=D_1+X_0=\pi(D_1)+X_0$, $D$ is a convex polyhedron in $W$ by Proposition~\ref{Zheng_ext}. Hence there exist $w^*_1, \dots, w^*_m \in W^*$ and $\alpha_1, \dots, \alpha_m \in \mathbb R$ such that 
\begin{equation*}
D=\left\{x \in W \mid \; \langle w^*_i, x \rangle \leq \alpha_i,\ i=1,\dots,m \right\}.
\end{equation*}  
According to Lemma \ref{Hann_Banach}, there exist $x^*_i \in X^*$, $i=1,\dots,m$, such that $\langle x^*_i, x \rangle = \langle w^*_i, x \rangle$ for all $x \in W.$  
Therefore
\begin{equation*}
D=\left\{x \in W \mid \; \langle x^*_i, x \rangle \leq \alpha_i,\ i=1,\dots,m \right\}.
\end{equation*}
It follows that $D$ is a generalized polyhedral convex set in $X$.  $\hfill\Box$

\medskip
The next example is an illustration for Theorem \ref{main_result}. 
\begin{Example}\label{Ex1}{\rm
Let $X=C[a, b]$ be the linear space of continuous real valued functions on the interval $[a, b]$ with the norm defined by $||x||=\max_{t \in [a, b]} |x(t)|$. By the Riesz representation theorem (see e.g. \cite[Theorem~6, p.~374]{Kolmogorov_Fomin_1975} and \cite[Theorem~1, p.~113]{Luenberger_1969}), the dual space of $X$ is the {\it normalized space $X^*=NBV[a, b]$ of functions of bounded variation}, that is functions $y: [a, b] \rightarrow \mathbb{R}$ of bounded variation, $y(a)=0$, and $y(\cdot)$ is continuous from the left at every point of $(a, b)$.  Let $x^*_1, x^*_2 \in X^*$ be defined by
\begin{equation}\label{R_integral}
\langle x^*_1, x \rangle =\int\limits_a^b \omega_1(t) x(t)dt,   \quad \langle x^*_2, x \rangle =\int\limits_a^b \omega_2(t) x(t)dt,
\end{equation}
where $\omega_1, \omega_2$ in $X\setminus\{0\}$ are chosen such that the vectors $\omega_1, \omega_2$ are linearly independent. The integrals in \eqref{R_integral} are Riemannian. They equal respectively to the Riemann-Stieltjes integrals (see \cite[p.~367]{Kolmogorov_Fomin_1975})
$\int\limits_a^b x(t)dy_1(t)$ and $\int\limits_a^b x(t)dy_2(t)$, which are 
given by the $C^1$-smooth functions $y_i(t)=\int\limits_a^t \omega_i(\tau)d \tau$, $i=1,2$. 
   Set
\begin{equation*}
\gamma_{i,j}:=\int\limits_a^b \omega_i(t) \omega_j(t)dt, 
\end{equation*}   
for $i, j \in \{1,2\}$.
It is clear that $\gamma_{1,2}=\gamma_{2,1}, \gamma_{1,1} >0, \gamma_{2,2} >0$. The Cauchy-Schwarz inequality
$$\left(\int\limits_{a}^{b} x^2(t)dt \right)^{1/2}\left(\int\limits_{a}^{b} y^2(t)dt \right)^{1/2} \geq \left|\int\limits_{a}^{b} x(t)y(t)dt \right|,$$
which is valid for any functions $x(\cdot), y(\cdot) \in C[a, b] \subset L_2[a, b]$, implies that $\delta:=\gamma_{1,1}.\gamma_{2,2} - \gamma_{1,2}^2 \geq 0$. As the vectors $\omega_1, \omega_2$ are linearly independent, we must have $\delta >0$.
Given any $\alpha_1, \alpha_2 \in \mathbb{R}$, we want to find a representation of form \eqref{rep_2} for the convex polyhedron
\begin{equation}\label{formula_for_D_Ex1}
D:=\left\{x \in X \mid \langle x^*_1, x \rangle \leq \alpha_1, \; \langle x^*_2, x \rangle \leq \alpha_2\right\}.
\end{equation}
Let $X_0:=\left\{x \in X \mid \langle x^*_1, x \rangle =0, \; \langle x^*_2, x \rangle =0 \right\}.$
For $x=\eta_1 \omega_1 + \eta_2 \omega_2$ with $\eta_1, \eta_2 \in \mathbb{R}$, we have
\begin{equation*}
	\langle x_1^*, x \rangle = \eta_1 \int\limits_{a}^{b} \omega_1^2 (t)dt + \eta_2 \int\limits_{a}^{b} \omega_1(t) \omega_2(t) dt=\eta_1 \gamma_{1,1} + \eta_2 \gamma_{1,2},
\end{equation*}
and
\begin{equation*}
\langle x_2^*, x \rangle = \eta_1 \int\limits_{a}^{b} \omega_1(t)\omega_2(t)dt + \eta_2 \int\limits_{a}^{b} \omega_2^2(t) dt=\eta_1 \gamma_{1,2} + \eta_2 \gamma_{2,2}.
\end{equation*}
Since $\delta > 0$, there exists an unique pair of real numbers $(\eta_1, \eta_2)$ satisfying 
\begin{equation*}
\begin{cases}
\eta_1 \gamma_{1,1} + \eta_2 \gamma_{1,2}=\alpha_1\\
\eta_1 \gamma_{1,2} + \eta_2 \gamma_{2,2}=\alpha_2.
\end{cases}
\end{equation*}
Let the point $u$ and the directions $v_1, v_2 \in X$ be defined by
\begin{equation*}
u=\eta_1 \omega_1+\eta_2 \omega_2, \quad v_1=\gamma_{1,2}\omega_1 - \gamma_{1,1}\omega_2, \quad v_2=-\gamma_{2,2}\omega_1+\gamma_{1,2}\omega_2.
\end{equation*}
It is easy to verify that $\langle x^*_i, u \rangle = \alpha_i$ for $i=1,2$, and
\begin{equation*}
\langle x^*_1, v_1 \rangle = 0, \quad \langle x^*_1, v_2 \rangle = -\delta, \quad \langle x^*_2, v_1 \rangle = -\delta, \quad \langle x^*_2, v_2 \rangle = 0.
\end{equation*}
Let us show that 
\begin{equation}\label{D_Ex1}
D=\{u+\mu_1 v_1 + \mu_2 v_2 \mid \mu_j \geq 0,\, j=1,2\}+X_0.
\end{equation}
Take any $x=u+\mu_1 v_1 + \mu_2 v_2 +x_0$ with $\mu_1, \mu_2 \in \mathbb R_+$ and $x_0 \in X_0$. Because
\begin{equation*}
\langle x_1^*, x \rangle = \langle x_1^*,u \rangle + \mu_1 \langle x_1^*,v_1 \rangle+\mu_2 \langle x_1^*,v_2 \rangle+\langle x_1^*,x_0 \rangle=\alpha_1-\mu_2 \delta \leq \alpha_1
\end{equation*}
and
\begin{equation*}
\langle x_2^*, x \rangle= \langle x_2^*,u \rangle + \mu_1 \langle x_2^*,v_1 \rangle+\mu_2 \langle x_2^*,v_2 \rangle+\langle x_2^*,x_0 \rangle=\alpha_2-\mu_1 \delta \leq \alpha_2,
\end{equation*}
we have $x \in D$. Now, take any $x \in D$. Put $\mu_1=\delta^{-1}\left( \alpha_2 - \langle x_2^*, x \rangle\right)$, $\mu_2=\delta^{-1}\left( \alpha_1 - \langle x_1^*, x \rangle\right)$, and $x_0=x-\left( u+\mu_1 v_1 + \mu_2 v_2\right) $. Note that $\mu_1 \geq 0$, $\mu_2 \geq 0$ and $x=u+\mu_1 v_1 + \mu_2 v_2+x_0$. Since $\langle x^*_i, x_0 \rangle =0$ for $i=1,2$, we see that $x_0 \in X_0$. The formula \eqref{D_Ex1} has been proved.}
\end{Example}

Based on the preceding example, we can easily construct an illustrative example for polyhedral convex sets in locally convex Hausdorff topological vector spaces.

\begin{Example}\label{Ex2} {\rm
Keeping all the notations of Example \ref{Ex1}, we consider $X=C[a, b]$ with the weak topology. Then $X$ is a locally convex Hausdorff topological vector space whose topology is not a norm  topology. The analysis given above shows that the set $D$ in \eqref{formula_for_D_Ex1} admits the representation \eqref{rep_2}.	}
\end{Example}

From Theorem \ref{main_result} we can obtain a representation formula for generalized polyhedral convex cones.
\begin{Theorem}\label{gpc_cone_representation} A nonempty set $K \subset X$ is a generalized polyhedral convex cone if and only if there exist $v_j \in K, j=1,\dots, \ell$, and a closed linear subspace $X_0$ such that     
\begin{equation}\label{rep_3}
K=\Bigg\{\sum\limits_{j=1}^\ell \mu_j v_j \mid  \mu_j \geq 0, \ \forall j=1,\dots,\ell   \Bigg\}+X_0. 
\end{equation}
\end{Theorem}
{\bf Proof.} {\it Necessity:} If $K$ is a generalized polyhedral convex cone, then by Theorem \ref{main_result} we can find $u_i \in K$, $i=1,\dots,k$, $v_j \in X$, $j=1,\dots,\ell$, 
and a closed linear subspace $X_0$ such that
\begin{equation}\label{rep_2a}
\begin{aligned}
K=\Bigg\{ \sum\limits_{i=1}^k \lambda_i u_i + \sum\limits_{j=1}^\ell \mu_j v_j \mid & \lambda_i \geq 0, \ \forall i=1,\dots,k,    \\ 
&\sum\limits_{i=1}^k \lambda_i=1,\ \, \mu_j \geq 0,\ \forall j=1,\dots,\ell   \Bigg\}+X_0.&& 
\end{aligned}
\end{equation}
To show that $v_j \in K$ for $j=1,\dots,\ell$, it suffices to observe by \eqref{rep_2a} that $\frac{1}{t}\left(u_1+tv_j\right)=\frac{1}{t}u_1+v_j$ belongs to $K$ for all $t > 0$, because $K$ is a cone. Letting $t \to \infty$, by the closedness of $K$, we get $v_j \in K$. Since $v_j \in K$ for $j=1,\dots,\ell$, and since $t_iu_i \in K$ for all  $i=1,\dots,k$, and $t_i \geq 0$, by choosing $v_{\ell+i}=u_i$ for $i=1,\dots,k,$        
by \eqref{rep_2a} we see  that $K$ admits the representation \eqref{rep_3} where $\ell$ is replaced by $\ell+k$.

{\it Sufficiency:} If $K$ has the form \eqref{rep_3} then it is a cone. In addition, $K$ is a generalized polyhedral convex set by Theorem \ref{main_result}.$\hfill\Box$

\section{Solution Existence in Linear Optimization}
Consider a {\it generalized linear programming problem} 
\begin{equation*}
{\rm (LP)} \qquad \min \left\{\langle x^*, x \rangle \mid x \in D\right\} 
\end{equation*}
where, as before, $X$ is a locally convex Hausdorff topological vector space, $D \subset X$ is a generalized polyhedral convex set, $x^* \in X^*$. By definition, the {\it recession cone} $0^{+}C$ of a convex set $C \subset X$ is given by 
\begin{equation}
0^{+}C=\{v \in X \mid x+tv \in C, \ \forall x \in X,\ \forall t \in \mathbb{R}_+\}.
\end{equation} 
If $D$ is represented in the form \eqref{rep_2}, then $0^+D={\rm{cone}}\{v_1, \dots,v_{\ell}\}+X_0$.

\medskip 
The following two existence theorems for {\rm (LP)} are known results. Actually, in combination, they express the contents of Theorem 2.199 from \cite{Bonnans_Shapiro_2000}. The latter, in its turn, is a special case of Theorem 2.198 from \cite{Bonnans_Shapiro_2000}. The simple proofs given below show how Theorem \ref{main_result} can be used to study the solution existence of generalized linear programs.     

\begin{Theorem}\label{Eaves} {\rm (The Eaves-type Existence Theorem; see \cite[Theorem 2.199]{Bonnans_Shapiro_2000})} If $D$ is nonempty, then {\rm(LP)} has a solution if and only if $\langle x^*, v \rangle \geq 0$ for every $v \in 0^+D$.
\end{Theorem}
{\bf Proof.} If {\rm(LP)} has a solution $\overline{x}$, then for each $v \in 0^+D$ it holds that
\begin{equation*}
\langle x^*, \overline{x} \rangle \leq \langle x^*, \overline{x}+tv \rangle=\langle x^*, \overline{x}\rangle + t  \langle x^*, v \rangle \ \; \forall t \in \mathbb{R}_+, 
\end{equation*}
because $\overline{x}+tv \in D$ for every $t \geq 0$. Hence $\langle x^*, v \rangle \geq 0$.

Conversely, suppose that	$\langle x^*, v \rangle \geq 0$ for every $v \in 0^+D$. Let us represent $D$ in the form \eqref{rep_2}. Select an element $u_{i_0} \in \{u_1, \dots, u_k\}$ such that 
\begin{equation*}
\langle x^*, u_{i_0} \rangle =\min\Big\{ \langle x^*, u_i \rangle \mid i=1,\dots,k\Big\}.
\end{equation*} 
By \eqref{rep_2}, for every $x \in D$ there exist $\lambda_i \in \mathbb{R}_+$, $i=1,\dots,k$, $\sum\limits_{i=1}^k \lambda_i=1$, and $v \in 0^+D$ such that $x=\sum\limits_{i=1}^k\lambda_i u_i + v$. Then we have
\begin{equation*}
\langle x^*, x \rangle= \sum\limits_{i=1}^k \lambda_i \langle x^*,u_i \rangle + \langle x^*,v \rangle 
\geq \sum\limits_{i=1}^k \lambda_i \langle x^*,u_{i_0} \rangle = \langle x^*,u_{i_0} \rangle.
\end{equation*}  
Since $x \in D$ can be chosen arbitrarily, $u_{i_0}$ must be a solution of {\rm (LP)}. $\hfill\Box$

\begin{Remark}{\rm
If $\langle x^*, v \rangle \geq 0$ for every $v \in 0^+D$, then one says that the functional $x^*$ is {\it copositive} on the recession $0^+D$. We called Theorem \ref{Eaves} the Eaves-type Existence Theorem in linear optimization to trace back Eaves' idea \cite[Theorem 3 and Corollary 4, p.702]{Eaves_1971} (see also \cite[Theorem 2.2]{Lee_Tam_Yen_2005}) in using recession cones for existence theorems in quadratic programming.       
}\end{Remark}
 
\begin{Theorem}\label{Frank_Wolfe} {\rm (The Frank--Wolfe-type Existence Theorem; see \cite[Theorem 2.199]{Bonnans_Shapiro_2000})} If $D$ is nonempty, then {\rm(LP)} has a solution if and only if there is a real number $\gamma$ such that $\langle x^*, x \rangle \geq \gamma$ for every $x \in D$.
\end{Theorem}
{\bf Proof.} The necessity  is obvious. To prove the sufficiency, suppose that there is a $\gamma \in \mathbb{R}$ such that $\langle x^*, x \rangle \geq \gamma$ for all $x \in D$. Then, for any $v \in 0^+D$ and $x \in D$ we have
\begin{equation*}
\gamma \leq \langle x^*, x+tv \rangle=\langle x^*, x \rangle+t\langle x^*, v \rangle
\end{equation*}
 for every $t > 0$. It follows that $\langle x^*, v \rangle \geq 0$ for any $v \in 0^+D$. So,  by Theorem \ref{Eaves}, we can assert that {\rm(LP)} has a solution.
$\hfill\Box$

\begin{Remark}{\rm Due to the formulation of the existence theorem in quadratic programming of Frank and Wolfe \cite[p. 158]{Frank_Wolfe_1956} (see also \cite[Theorem 2.1]{Lee_Tam_Yen_2005}), we called Theorem~\ref{Frank_Wolfe} the Frank--Wolfe-type Existence Theorem in linear optimization.
	}\end{Remark}

We are interested in studying the region $G$ of all $x^*$ for which {\rm (LP)} has a nonempty solution set, assuming that the constraint set $D$ is nonempty and fixed. 
\begin{Proposition}\label{description_G}
If $D$ has the form \eqref{rep_2}, then $G$ is a generalized polyhedral convex cone of $X^*$ which has the representation
\begin{equation}\label{reps_domain_G}
G=X_0^{\perp} \cap \{x^* \in X^* \mid \langle x^*, v_j \rangle \geq 0,\  j=1,\dots,\ell \}
\end{equation}
where  $X_0^{\perp}:=\{x^* \in X^* \mid \langle x^*, x_0 \rangle = 0, \ \forall x_0 \in X_0\}$ is the annihilator {\rm \cite[p.~117]{Luenberger_1969}} of $X_0$. 
\end{Proposition}
{\bf Proof.} By Theorem \ref{Eaves}, $G=\{x^* \in X^* \mid \langle x^*, v \rangle \geq 0, \ \forall v \in 0^+D\}$. Therefore, given any $x^* \in G$, we have $\langle x^*, v \rangle \geq 0$ for all $v \in 0^+D$. Hence, for every $x_0 \in X_0 \subset 0^+D$ one has $\langle x^*, x_0 \rangle =0$ because $-x_0 \in 0^+D$. Thus $x^* \in X_0^{\perp}$. In addition, for each $j=1,\dots,\ell$, one has $\langle x^*, v_j \rangle \geq 0$ as $v_j \in 0^+D$. This establishes the inclusion ``$\subset$'' in \eqref{reps_domain_G}. 

Conversely, suppose that $x^* \in X_0^{\perp} \cap \{x^* \in X^* \mid \langle x^*, v_j \rangle \geq 0,\ j=1,\dots,\ell \}$. Since $0^+D={\rm{cone}}\{v_1, \dots,v_{\ell}\}+X_0$, the last inclusion implies that $\langle x^*, v \rangle \geq 0$, for all $v \in 0^+D$. Hence, by Theorem \ref{Eaves} we can conclude that $x^* \in G$. The inclusion ``$\supset$'' in \eqref{reps_domain_G} has been proved.

From \eqref{reps_domain_G} it follows that $G$ is a generalized polyhedral convex set. The fact that $G$ is a cone is obvious.  $\hfill\Box$

\medskip
Next, for each $x^* \in G$, we want to describe the solution set of {\rm (LP)}, which is denoted by $S(x^*)$. For doing so, let us suppose that $D$ is given by \eqref{rep_2} and consider the index sets 
$$I(x^*):=\{i_0 \in \{1,\dots,k\} \mid \langle x^*, u_{i_0} \rangle \leq \langle x^*, u_{i} \rangle \ \, \forall i=1,\dots,k \},$$
and $$J(x^*):=\{j_0 \in \{1,\dots,\ell\} \mid \langle x^*, v_{j_0} \rangle =0 \}.$$
Note that $I(x^*)$ is nonempty, but it may happen that $J(x^*)$ is empty.  
\begin{Proposition}	
	If $x^* \in G$ and $D$ is given by \eqref{rep_2}, then 
	\begin{equation}\label{solution_map_S}
	\begin{aligned}
	S(x^*)=\Bigg\{ \sum\limits_{i \in I(x^*)} \lambda_i u_i &+ \sum\limits_{j \in J(x^*)} \mu_j v_j \mid \ \lambda_i \geq 0 \ \, \forall i \in I(x^*),    \\ 
	&\sum\limits_{i \in I(x^*)} \lambda_i=1,\; \mu_j \geq 0 \ \, \forall j \in J(x^*)   \Bigg\}+X_0.&& 
	\end{aligned}
	\end{equation}
	In particular, $S(x^*)$ is a generalized polyhedral convex set. 
\end{Proposition}
{\bf Proof.} First, take an arbitrary element $\bar{x}$ from the set on the right-hand-side of \eqref{solution_map_S}. Let 
$$\bar{x}=\sum\limits_{i \in I(x^*)} \bar{\lambda}_i u_i + \sum\limits_{j \in J(x^*)} \bar\mu_j v_j+\bar x_0,$$
where $\bar\lambda_i \geq 0$ for all $i \in I(x^*)$, $\sum\limits_{i \in I(x^*)} \bar\lambda_i=1$, $\bar\mu_j \geq 0$ for all $j \in J(x^*)$, and $\bar x_0 \in X_0$. By \eqref{rep_2}, for each $x \in D$ one can find $\lambda_i \geq 0$ for $i=1,\dots,k$, $\sum\limits_{i=1}^k \lambda_i=1$, $\mu_j \geq 0$ for $j=1,\dots,\ell$, and $x_0 \in X_0$ such that
 $$x=\sum\limits_{i=1}^k {\lambda}_i u_i + \sum\limits_{j=1}^{\ell} \mu_j v_j+x_0.$$
By Proposition \ref{description_G}, $\langle x^*, x_0 \rangle = \langle x^*, \bar x_0 \rangle =0$.
If $J(x^*) \neq \emptyset$, then using Theorem \ref{Eaves} and formula $0^+D={\rm{cone}}\{v_1, \dots,v_{\ell}\}+X_0$ we get 
$$\left\langle x^*, \sum\limits_{j=1}^{\ell} \mu_j v_j \right\rangle \geq 0=  \left\langle x^*, \sum\limits_{j \in J(x^*)}\bar\mu_j v_j \right \rangle. $$
Now, selecting an index ${i_0} \in I(x^*)$ and recalling the definition $I(x^*)$, we get
\begin{equation*}
\begin{aligned}
\left\langle x^*, \sum\limits_{i=1}^{k} \lambda_i u_i \right\rangle \geq \sum\limits_{i=1}^{k} \lambda_i \left\langle x^*, u_{i_0} \right\rangle & =\left\langle x^*, u_{i_0} \right\rangle  \\
&=\left\langle x^*,\sum\limits_{i \in I(x^*)} \bar\lambda_iu_i \right\rangle.
\end{aligned}
\end{equation*}
It follows that $\langle x^*, x \rangle \geq \langle x^*, \bar{x} \rangle$. We have shown that $\bar{x} \in S(x^*)$.

Second, take any vector $\bar x \in S(x^*)$ and represent it in the form
$$\bar{x}=\sum\limits_{i=1}^k \bar{\lambda}_i u_i + \sum\limits_{j=1}^{\ell} \bar\mu_j v_j+\bar x_0,$$
where $\bar\lambda_i \geq 0$ for $i=1,\dots,k$, $\sum\limits_{i=1}^k \bar\lambda_i=1$, $\bar\mu_j \geq 0$ for $j=1,\dots,\ell$, and $\bar x_0 \in X_0$. It is easy to show that $\bar\lambda_i=0$ for all $i \notin I(x^*)$ and $\bar\mu_j=0$ for all $j \notin J(x^*)$. This implies that $\bar x$ belongs to the set on the right-hand-side of \eqref{solution_map_S}.  

The proof is complete. $\hfill\Box$

\section{The Weakly Efficient Solution Set in Linear Vector Optimization}
Consider a {\it generalized linear vector optimization problem} of the form  
\begin{equation*}
{\rm (VLP)} \qquad \quad {{\rm min}_K} \left\{Mx \mid x \in D\right\}
\end{equation*}
with $M: X \rightarrow Y$ being a continuous linear mapping between locally convex Hausdorff topological vector spaces, $D \subset X$ a generalized polyhedron, $K \subset Y$ a polyhedral convex cone. 

\medskip
We say that $u \in D$ is a \textit{weakly efficient solution} of {\rm (VLP)} if there does not exist any $x \in D$ such that $Mu - Mx \in {\rm int}K$. The set of all the weakly efficient solutions is denoted by $E^w$. We are interested in finding conditions to have $E^w \neq \emptyset$.

\medskip
By a standard scalarization scheme in vector optimization, given any $y^* \in Y^*$, we define the scalar problem
\begin{equation*}
{\rm (LP)}_{y^*} \qquad \min \left\{\langle y^*, Mx \rangle \mid x \in D\right\}.
\end{equation*}

To make our presentation easier for reading, we give simple proof for the following known result.   
\begin{Lemma}\label{Luc} {\rm (See \cite[Proposition~3.2, p.~95]{Luc})} If ${\rm int}K$ is nonempty, then $u \in D$ is a weakly efficient solution of {\rm (VLP)} if and only if there exists $y^* \in K^* \setminus\{0\}$, where $K^*:=\left\lbrace y^* \in Y^* \mid \langle y^*, y \rangle \geq 0 \ \, \forall y \in K \right\rbrace$, such that 
\begin{equation}\label{argmin_scalarized}
u \in {\rm argmin} \left( {\rm (LP)}_{y^*}\right).
\end{equation}	
	
\end{Lemma}
{\bf Proof.} First, suppose that $u \in E^w$. Since $\big(Mu-M(D)\big) \cap {\rm int}K=~\emptyset$ and since $Mu-M(D)$ and ${\rm int}K$ are convex sets, by the separation theorem (see, e.g., \cite[Theorem 2.13]{Bonnans_Shapiro_2000}), there exists $y^* \in Y^* \setminus\{0\}$ such that
\begin{equation*}
\langle y^*, Mu-Mx \rangle \leq \langle y^*, y \rangle
\end{equation*}   
for all $x \in D$ and $y \in K$. Substituting $x=u$ to the above inequality yields $\langle y^*, y \rangle \geq 0$ for all $y \in K$. Hence $y^* \in K^*$. Choosing $y=0$, one has $\langle y^*, Mu \rangle \leq \langle y^*, Mx \rangle$ for every $x \in D$. This shows that the inclusion \eqref{argmin_scalarized} is valid.

Now, suppose that $u \in D$ and there is $y^* \in K^* \setminus\{0\}$ such that \eqref{argmin_scalarized} is satisfied. If $u \not\in E^w$, then there exist $x \in D$ and a balanced neighborhood $V$ of 0 satisfying $Mu-Mx + V  \subset K$. Hence, for each $v \in V$, one has $\langle y^*, Mu-Mx+v \rangle \geq 0$. In combination with the inequality $\langle y^*, Mu-Mx \rangle \leq 0$ which is guaranteed by \eqref{argmin_scalarized}, this implies $\langle y^*, v \rangle \geq 0$. As $V$ is a balanced neighborhood of 0, we can assert that $\langle y^*, v \rangle = 0$ for all $v \in V$. Let $y \in Y$ be such that $\langle y^*, y \rangle \not=0$. Since there exists $t >0$ with $ty \in V$, we get $\langle y^*, y \rangle=0$, a contradiction.      
$\hfill\Box$

\begin{Remark}{\rm Looking back to the proof of Lemma \ref{Luc}, we can observe that it suffices to assume that $K \subset Y$ is a convex cone. In other words, the polyhedrality of $K$ is superfluous for the assertion of the lemma.} 
\end{Remark}

We have $\langle y^*, Mx \rangle = \langle M^*y^*, x \rangle$ for all $x \in X$, where $M^*: Y^* \rightarrow X^*$ is the adjoint operator of $M$. 

\begin{Theorem}\label{VLP_existence_theorem_1}
  Problem {\rm (VLP)} has a solution if and only if 
  \begin{equation}\label{criterion_1} 
  M^*\left(K^*\setminus\{0\}\right) \cap \left(0^+D\right)^* \neq \emptyset.
  \end{equation}
  In particular, if 
  \begin{equation}\label{criterion_2} 
  M^*\left(K^*\right) \cap \left(0^+D\right)^* \neq \{0\},
  \end{equation}
   then {\rm (VLP)} has a solution.
\end{Theorem}
{\bf Proof.} By Lemma \ref{Luc}, {\rm (VLP)} has a solution if and only if there exists $y^* \in K^* \setminus\{0\}$ such that the solution set of ${\rm (LP)}_{y^*}$ is non empty. According to Theorem \ref{Eaves}, this solution set is non-void if and only if $M^*y^* \in \left(0^+D\right)^*$. Thus, we have shown that {\rm (VLP)} has a solution if and only if \eqref{criterion_1} is fulfilled.     

Now, suppose that \eqref{criterion_2} is satisfied. Then we can find $y^* \in K^*$ such that $M^*y^* \in \left(0^+D\right)^*$ and $M^*y^* \neq 0$. Since the later obviously implies that $y^* \neq 0$, we have $y^* \in M^*\left(K^*\setminus\{0\}\right) \cap \left(0^+D\right)^*$. Hence \eqref{criterion_1} holds true, so {\rm (VLP)} has a solution.~$\hfill\Box$      

\begin{Remark}{\rm The assertions of Theorem \ref{VLP_existence_theorem_1} are valid for the case $K \subset Y$ is an arbitrary convex cone.} 
\end{Remark}
  
We conclude this section by a statement about the structure of $E^w$ which is applicable also the case $K \subset Y$ is an arbitrary convex cone.

\begin{Theorem}\label{structure_solution_set_1} The weakly efficient solution set $E^w$ of {\rm (VLP)} is the union of finitely many generalized polyhedral convex sets.
\end{Theorem}
{\bf Proof.} Using Lemma \ref{Luc}, we can represent the weakly efficient solution set of {\rm (VLP)} as follows 
\begin{equation}\label{repr_solution_set_1}
E^w=\bigcup\limits_{y^* \in K^*\setminus\{0\}} {\rm argmin} \left( {\rm (LP)}_{y^*}\right).
\end{equation}
Setting $x^*=M^*y^*$, we can rewrite \eqref{repr_solution_set_1} as
\begin{equation}\label{repr_solution_set_1a}
E^w=\bigcup\limits_{x^* \in M^*\left( K^*\setminus\{0\}\right) } S(x^*),
\end{equation}
where $S(x^*)$ is the solution set of the problem {\rm (LP)} considered in Section 3. Invoking \eqref{solution_map_S} and noting that the number of the index sets $I(x^*)$ (resp., the number of the index sets $J(x^*)$) is finite, from \eqref{repr_solution_set_1a} we obtain the desired conclusion. $\hfill\Box$

\end{document}